\documentclass[12pt]{article}
\usepackage{theorem, amsfonts,amsmath}
\newtheorem{lemma}{Lemma}[section]

\newtheorem{theorem}[lemma]{Theorem}
\newtheorem{corollary}[lemma]{Corollary}
{\theorembodyfont{\upshape}\newtheorem{definition}[lemma]{Definition}}
{\theorembodyfont{\upshape}
{\theorembodyfont{\upshape}\newtheorem{remark}[lemma]{Remark}}
{\theorembodyfont{\upshape}}
{\theorembodyfont{\upshape}}

\newcommand{\ep}{\varepsilon}
\newcommand{\Proof}{\underbar{Proof}{\hskip 0.1in}}
\newcommand{\rn}{\mathbb{R}^N}
\newcommand{\hf}{\frac{1}{2}}
\newcommand{\qt}{\frac{1}{4}}

\newcommand{\id}{{\rm d}}
\newcommand{\dvol}{{\rm d}vol}
\newcommand{\dnx}{{\rm d}^Nx}
\newcommand{\dn}{{\rm d}^N}
\newcommand{\bd}{\partial \Omega}

\begin{document}
\title{LOG--SOBOLEV INEQUALITIES AND REGIONS WITH EXTERIOR EXPONENTIAL CUSPS}
\author{C. Mason}
\date{June 2001}
\maketitle
\begin{abstract}  We begin by studying certain semigroup estimates which are more singular than those implied by a Sobolev embedding theorem but which are equivalent to certain logarithmic Sobolev inequalities.  We then give a method for proving that such log--Sobolev inequalities hold for Euclidean regions which satisfy a particular Hardy--type inequality.  Our main application is to show that domains which have exterior exponential cusps, and hence have no Sobolev embedding theorem, satisfy such heat kernel bounds provided the cusps are not too sharp.  Finally we consider a rotationally invariant domain with an exponentially sharp cusp and prove that ultracontractivity breaks down when the cusp becomes too sharp.
\end{abstract}

\section{Introduction}

The spectral behaviour of the Neumann Laplacian, $H_N$, is known to be extremely sensitive to the regularity of the boundary.  There is a substantial body of research that shows how to produce peculiar behaviour.  We mention in particular the work of Simon and his various co--authors \cite{S92,HSB} and also Evans and Harris \cite{EH87} (futher references can be found in these papers).

In the opposite direction the spectrum can be shown to be well behaved if one can show that the associated semigroup $e^{-H_Nt}$  is {\em ultracontractive}; i.e. it is bounded from $L^2$ to $L^\infty$ for $0<t\leq 1$.  In the case that the space has finite measure this implies for example that the resolvent is compact, that the associated eigenfunctions all lie in $L^\infty$ and that $e^{-H_Nt}$ is compact on $L^p$ for all $1\leq p\leq \infty$ and $0<t<1$.

A further motivation for proving such results is the following.  No matter how bad the theoretical results can be, numerical methods of computing Neumann eigenvalues must assume that if one region is approximated by another then eigenvalues will still exist and will be close to those of the original region.  Burkenov and Davies, \cite{SSNL}, consider this problem and are able to give precise theorems which justify such methods for domains with H\"older class boundaries.  To do this they study the associated semigroup $e^{-H_Nt}$ and show that it is  ultracontractive.

In many cases (such as those considered by Burenkov and Davies) proving ultracontractivity  can be achieved by proving a Sobolev embedding of the form
\[
W^{1,2} \hookrightarrow L^q  
\]
for some $q>2$.  While this is often possible it fails in the case that the region has exterior exponentially sharp cusps (\cite[Theorem 5.32]{A}).

Of course, the lack of a Sobolev embedding theorem says nothing about the possibility of proving semigroup and heat kernel bounds and hence results about the spectrum.

In this paper we are motivated by the results of Davies and Burenkov to study the question of how singular a domain can be and still possess an ultracontractive estimate.  Moreover, we will investigate the implications these results have for bounds on the eigenvalues and eigenfunctions and also the use of Hardy--type inequalities.  However, we will  use a more general tool than the Sobolev embedding, namely the {\em logarithmic Sobolev inequality}.  We begin in section \ref{sect:logsobinequalities}  by studying the type of inequalities that will be proved and some  spectral consequences;  Theorem \ref{lem:eigenestimates} gives lower bounds on the rate at which the eigenvalues grow and upper bounds on the $L^\infty$ norm of the eigenfunctions.  Conversely we show that we show that these bounds imply a log--Sobolev inequality -- see Theorem \ref{thm:kernelestimate} and Corollary \ref{cor:semigroupestimate} .

Our main tool in actually proving the inequalities is to first prove in Theorem \ref{thm:gls} a generalised log--Sobolev inequality that is valid for arbitrary bounded regions in $\rn$.  This will then be combined with a Hardy--type inequality which we study in section \ref{sect:loghardy}.  As an example we consider in section \ref{sect:euclideandomains} a simple region that may have exterior exponential cusps and show that the associated Neumann semigroup is indeed ultracontractive provided the cusp is not too sharp.  

Our final result, in section \ref{sect:rot} is to consider a rotationally invariant domain with an exponentially sharp cusp which shows that ultracontractivity {\em does} break down if the cusp is too sharp.  

\section{Log--Sobolev Inequalities}\label{sect:logsobinequalities}

Let $\Omega$ be a region  in $\rn$ and define the Neumann Laplacian to be the non--negative self--adjoint operator $H_N$ acting in $L^2(\Omega)$ associated with the quadratic form
\[
Q(f)=\left\{\begin{array}{cc} \int_\Omega |\nabla f|^2\dnx & \text{ if } f \in W^{1,2}(\Omega) \\ +\infty & \text{otherwise.}\end{array} \right.
\]
The associated symmetric Markov semigroup is denoted by $e^{-H_Nt}$.

We refer to \cite[Chapter 2]{HKST} for an introduction to logarithmic Sobolev inequalities.  The following theorem captures one of the main results, namely that a log--Sobolev inequality with a suitable right hand side is equivalent to an ultracontractive estimate.

\begin{theorem}\label{thm:equivultra}
Let $\alpha>1$.  Then the following are equivalent.
\begin{enumerate}
\item The log--Sobolev inequality
\[
\int_\Omega f^2\log_+ f \leq \ep Q(f)+\eta(\ep)\|f\|_2^2+\|f\|_2^2\log \|f\|_2
\]
is valid for $0\leq f\in W^{1,2}(\Omega), 0<\ep<1$ and there exists $c_1>0$ such that $\eta$ satisfies 
\[
\eta(\ep) \leq c_1\ep^{-1/(\alpha-1)}.
\]
\item The  semigroup $e^{-H_Nt}$ satisfies
\begin{align}\label{eqn:normsemi}
\|e^{-H_Nt}f\|_\infty \leq c_2\exp(c_3t^{-1/(\alpha-1)})\|f\|_2
\end{align}
for all $f\in L^2(\Omega)$, some constants $c_2,c_3>0$ and $0<t\leq 1$.
\item $e^{-H_Nt}$ has a continuous integral kernel $K(t,x,y)$ and there exist  $c_4,c_5>0$ such that
\begin{align}\label{eqn:kernelbound}
0< K(t,x,y) \leq c_4\exp(c_5t^{-1/(\alpha-1)})
\end{align}
for $0<t\leq 1$ and $x, y \in \Omega$.
\end{enumerate}
 \end{theorem}

We now turn our attention to the the eigenvalues and eigenfunctions.  Suppose $H_N$ has compact resolvent (which it has if $|\Omega|<\infty$ and any of the statements in Theorem \ref{thm:equivultra} hold) and denote its eigenvalues by
\[
0\leq \lambda_0 \leq \lambda_1 \leq \cdots \leq \lambda_n  \to \infty
\]
where we repeat each eigenvalue according to its multiplicity.  The associated orthonormal eigenfunctions are denoted by $f_n$.  The following lemma is an easy consequence of ultracontractivity.

\begin{lemma}\label{lem:eigenestimates} Let $\Omega$ be a region of finite measure such that any one of the statements in Theorem \ref{thm:equivultra} holds.  Then there exist $c_6,c_7>0$ and $N>c_7^{-1}$ such that 
\begin{align}\label{eqn:valuelower}
\lambda_n \geq c_6(\log (c_7n))^{\alpha}
\end{align}
for all $n\geq N$.  Also there exists $c_8,c_9>0$ such that
\begin{align}\label{eqn:vectorupper}
\|f_n\|_\infty \leq c_8 \left\{ \begin{array}{cc} 1 & 0\leq n < N \\ \exp (c_9\lambda_n^{1/\alpha}) & n \geq N \end{array} \right.
\end{align}
for all $n\geq N$.
\end{lemma}

\Proof  By integrating (\ref{eqn:kernelbound}) where $x=y$ over $\Omega$ we have
\[
ne^{-\lambda_nt} \leq \sum_{k=0}^n e^{-\lambda_kt} \leq \sum_{k=0}^\infty e^{-\lambda_kt}=\int_\Omega K(t,x,x)\dnx \leq c_4\exp(c_5t^{-1/(\alpha-1)})|\Omega|.
\]
Since $H_N$ has compact resolvent there exists $N$ such that $n\geq N$ implies that $\lambda_n\geq 1$.  For all such $\lambda_n$ put 
\[
t=\lambda_n^{-1+1/\alpha}
\]
 to get
\[
\frac{n}{c_4|\Omega|} \leq \exp ((c_5+1)\lambda^{1/\alpha}).
\]
This implies that 
\[
c_6(\log(nc_7))^{\alpha} \leq \lambda_n
\]
where $c_6=(c_5+1)^{-\alpha}$ and $c_7=(c_4  |\Omega|)^{-1}$ provided $nc_7\geq 1$ which we assume without loss of generality.

The second conclusion follows by putting $f=f_n$ into (\ref{eqn:normsemi}) for all $n\geq N$ to get
\[
e^{-\lambda_nt}\|f_n\|_\infty \leq c_2 \exp(c_3t^{-1/(\alpha-1)}).
\]
For $n \geq N$ set $t=\lambda_n^{-1+1/\alpha}$ to get
\[
\|f_n\|_\infty \leq c_2\exp ((c_3+1)\lambda^{1/\alpha}).
\]

For $n<N$ we put $t=1$ into  (\ref{eqn:normsemi}) to get
\[
\|f_n\|_\infty \leq c_2\exp(c_3)e.
\]
Thus we can take $c_8=c_2\exp(c_3+1)$ and $c_9=c_3+1$.

We now turn to the problem of proving a converse to the previous lemma.   First we give a simple but important lemma.

\begin{lemma}\label{lem:estbasic}
For $\alpha>1$ there exists $c_{10}>0$, depending only on $c_9$ and $\alpha>1$ such that
\[
\exp(-\lambda_nt/2+2c_9\lambda_n^{1/\alpha})\leq \exp (c_{10}t^{-1/(\alpha-1)}).
\]
\end{lemma}

\Proof  We use the inequality
\[
a\leq \ep a^\beta+\ep^{-1/(\beta-1)}
\]
valid for all $a,\ep>0$ and $\beta>1$, with $a=2c_9\lambda_n^{1/\alpha}$, $\ep c_9^{\alpha}2^{\alpha-1}=t$ and $\beta=\alpha$.  Thus
\[
2c_9\lambda_n^{1/\alpha}\leq \lambda_n t/2+2^{(\alpha+1)/(\alpha-1)}c_9^{\alpha/(\alpha-1)}  t^{-1/(\alpha-1)}
\]
from which the result follows.

\begin{theorem}\label{thm:kernelestimate} Suppose $H_N$ has discrete spectrum with non-negative eigenvalues $\lambda_n$ of finite multiplicity, written in increasing order and repeated according to multiplicity.  Let $f_n$ denote the corresponding orthonormal basis of eigenfunctions and suppose that inequalities (\ref{eqn:valuelower}) and (\ref{eqn:vectorupper}) are satisfied for  constants $c_6,c_7,c_8,c_9>0$, $N>c_7^{-1}$ and $\alpha>1$.  Then there exist constants $C_1,C_2>0$ depending only on  $c_6,c_7,c_8,c_9,N$ and $\alpha$ such that
\[
0<K(t,x,y)\leq C_1\exp(C_2t^{-1/(\alpha-1)})
\]
for $0<t\leq 1$.
\end{theorem}
\Proof  If $0<t\leq 1$ and $x,y\in \Omega$ then there exists $c_{11}>0$ such that
\begin{align}
0<K(t,x,y)&=\sum_{n=0}^\infty e^{-\lambda_nt}f_n(x)f_n(y) \notag\\
&\leq \sum_{n=0}^\infty e^{-\lambda_nt}\|f_n\|_\infty^2\notag \\
&\leq c_{8}\left(\sum_{n=0}^{N-1}1+\sum_{n=N}^\infty e^{-\lambda_nt+2c_9\lambda_n^{1/\alpha}}\right) \notag \\
&\leq  c_{8}\left(\sum_{n=0}^{N-1}1+\exp(c_{10}t^{-1/(\alpha-1)})\sum_{n=N}^\infty e^{-\lambda_nt/2}\right) \label{eqn:finalone},
\end{align}
where we apply Lemma \ref{lem:estbasic} to get the final line.  Next we observe that
\begin{align}
\sum_{n=N}^\infty e^{-\lambda_nt/2}&\leq \sum_{n=N}^\infty e^{-c_6(\log c_7n)^{\alpha}t/2} \notag \\
&\leq \int_N^\infty  e^{-c_6(\log c_7x)^{\alpha}t/2}\id x \notag \\
&=(c_7\alpha)^{-1}\int_{N'}^\infty e^{-c_6st/2+s^{1/(\alpha)}}s^{-(\alpha-1)/\alpha}\id s \notag \\
&\leq(c_7\alpha)^{-1}(\log c_7N)^{1-\alpha}\int_{N'}^\infty e^{-c_6st/2+s^{1/\alpha}} \id s \label{eqn:finaltwo}
\end{align}
where we have made the substitutions $s=(\log c_7x)^{\alpha}$ and $N'=(\log c_7 N)^{\alpha}$. 

A simple modification of Lemma \ref{lem:estbasic} gives us that
\[
s^{1/\alpha}\leq \frac{c_6ts}{4}+\left(\frac{c_6}{4}\right)^{-1/(\alpha-1)}t^{-1/(\alpha-1)}.
\]
Hence
\begin{align}
\int_{N'}^\infty e^{-c_6st/2+s^{1/\alpha}} \id s &\leq \exp \left(\left(\frac{c_6}{4}\right)^{-1/(\alpha-1)}t^{-1/(\alpha-1)}\right)\int_{N'}^\infty e^{-c_6st/4}\id s \notag\\
&=4\exp(c_{11}t^{-1/(\alpha-1)})\frac{e^{-c_6N't/4}}{c_6t} \notag \\
&\leq \frac{4}{c_6}\exp(c_{12}t^{-1/(\alpha-1)}) \label{eqn:finalfour},
\end{align}
where $c_{11}:=(c_6/4)^{-1(\alpha-2)}$ and $c_{12}>c_{11}+(\alpha-1)e^{-1}$.

Combining (\ref{eqn:finalone}),  (\ref{eqn:finaltwo}) and  (\ref{eqn:finalfour}) gives the final result.

 We have the following immediate corollary.

\begin{corollary}\label{cor:semigroupestimate}
Let the conditions of the previous theorem be satisfied.  Then there exist constants $C_3,C_4>0$ depending only on  $c_6,c_7,c_8,c_9,N$ and $\alpha$ such that
\[
\|e^{-H_Nt}f\|_\infty \leq C_3\exp\left(C_4t^{-1/(\alpha-1)}\right)
\]
for all $f\in L^2(\Omega)$ and $0<t\leq 1$.
\end{corollary}

\subsection{Generalised Log--Sobolev Inequality}

In this section we prove a generalised log--Sobolev inequality.  This will be our main tool in proving a log--Sobolev inequality for a region with exterior exponential cusps.  It is valid for arbitrary bounded regions in $\rn$.

\begin{definition}  Suppose $\bd\ne \emptyset$ and  define $d(x)$ to be the distance of $x$ from $\bd$.
\end{definition}

\begin{theorem}\label{thm:gls}
Let $\Omega$ be a domain in $\mathbb{R}^N$ with finite inradius.  For $0 \leq f \in W^{1,2}(\Omega)$ there exist constants $b_0,b_1,b_2>0$
\[
\int_\Omega f^2\log f \leq \ep Q(f)+\beta(\ep)\|f\|_2^2+\|f\|_2^2\log \|f\|_2^2+b_0\int_\Omega |\log d|f^2
\]
for all $\ep>0$ and some $\beta(\ep)$ satisfying
\[
\beta(\ep) \leq b_1-b_2 \log \ep.
\]
\end{theorem}

\Proof Given $\delta>0$ put
\[
S_\delta:=\{ x \in \Omega : 2\delta \leq d (x) \leq 3\delta \}
\]
and let $x_1, \dots ,x_{n(\delta)}$ be a maximal set of points in $S_\delta$ such that $|x_i-x_j|\geq \delta$ for all $i \ne j$.  For $\Omega$ bounded this number $n(\delta)$ is finite and moreover the number of balls containing  $x\in \Omega$ is bounded uniformly with respect to $x$ and $\delta$.

Now let $B_r(a)$ denote the ball centred at $a\in\Omega$ with radius $r>0$.  We then define the following norms and forms:
\begin{align*}
\| f\|_{2,a,r}^2&:=\int_{B_r(a)}|f|^2 \dnx \\
Q_{a,r}(f)&:=\int_{B_r(a)}|\nabla f|^2\dnx
\end{align*}
for $f\in W^{1,2}(B_{a,r})$.  Now given $a\in \mathbb{R}^N$ we have a log--Sobolev inequality for $B_1(a)$ namely 
\begin{align*}
\int_{B_1(a)}f^2\log_+ f \dnx &\leq \ep Q_{a,1}(f)+\tilde\beta(\ep)\|f\|_{2,a,1}^2+\|f\|_{2,a,1}\log \|f\|_{2,a,1}
\end{align*}
for $0 \leq f \in W^{1,2}(B_1(a))$
 and for all $\ep>0$ and $\tilde\beta(\ep)$ satisfying
\[
\tilde\beta(\ep)=b_3-\frac{N}{4}\log (\ep)
\]
for some constant $b_3>0$.  

By scaling we then have
\begin{align*}
\int_{B_\delta(a)}f^2\log_+ f\dnx &\leq \ep Q_{a,\delta}(f)+\tilde\beta(\ep)\|f\|_{2,a,\delta}^2+\|f\|_{2,a,\delta}\log \|f\|_{2,a,\delta}\\
&+(N/2)|\log \delta| \|f\|_{2,a,\delta}^2
\end{align*}

for $0\leq f \in W^{1,2}(B_\delta(a))$ and $\ep>0$.

Now suppose that $a\in S_\delta$.  For $x \in B_\delta(a)$ we have
\[
\delta \leq d(x) \leq 4 \delta.
\]
Thus if $\delta \geq 1$ then
\[
|\log d(x)| \geq |\log \delta| 
\]
and if $4\delta<1$ a calculation shows that
\[
(N/2+1)|\log d(x)| \geq (N/2)|\log \delta|.
\]
This range of $\delta$ will be sufficient for our purposes.

Thus if $\delta<1/4$ or $\delta \geq 1$ we have
\begin{align*}
\int_{B_\delta(a)}f^2\log_+ f \dnx &\leq \ep Q_{a,\delta}(f)+\tilde\beta(\ep)\|f\|_{2,a,\delta}^2+\|f\|_{2,a,\delta}\log \|f\|_{2,a,\delta} \\
&+(N/2+1)\int_{B_\delta(a)} |\log d||f|^2 \dnx.
\end{align*}

Given $\delta$ smaller than the inradius of $\Omega$ we have a natural restriction map
\[
R:W^{1,2}(\Omega)\to W^{1,2}(B_\delta(a))
\]
where
\[
(Rf)(x)=f(x).
\]
Hence given $0\leq f \in W^{1,2}(\Omega)$ and $\delta<1/4$ or $\delta \geq 1$ we have
\begin{align*}
\int_{B_\delta(a)}(Rf)^2\log_+ (Rf) \dnx &\leq \ep Q_{a,\delta}(Rf)+\tilde\beta(\ep)\|Rf\|_{2,a,\delta}^2+\|Rf\|_{2,a,\delta}\log \|Rf\|_{2,a,\delta} \\
&+(N/2+1)\int_{B_\delta(a)} |\log d||Rf|^2 \dnx.
\end{align*}

If we choose $\|f\|_{2}=1$ then $\|Rf\|_{2,a,\delta} \leq 1$ and hence
\begin{align*}
\int_{B_\delta(a)}(Rf)^2\log_+ (Rf) \dnx &\leq \ep Q_{a,\delta}(Rf)+\tilde\beta(\ep)\|Rf\|_{2,a,\delta}^2 \\
&+(N/2+1)\int_{B_\delta(a)} |\log d||Rf|^2 \dnx.
\end{align*}
We will now drop explicit reference to the restriction operator.

Let $\delta<1/4$ or $\delta \geq 1$ and $0\leq f\in W^{1,2}(\Omega)$ with $\|f\|_2=1$.  Then there exists $b_4>0$ 
\begin{align*}
\int_{S_\delta} f^2\log_+ f & \leq \sum_{i=1}^{n(\delta)} \int_{B_\delta(x_i)}f^2\log_+ f \\
& \leq b_4 \left(\int_{T_\delta}(\ep |\nabla f|^2+\tilde\beta(\ep)|f|^2 ) +(N/2+1)\int_{T_\delta} |\log d||f|^2 \right)
\end{align*}

 where $T_\delta:=S_{\delta/2}\cup S_{\delta}\cup S_{2\delta}$.  Now sum over $\delta=5^{-n}$ for all integers $n$ to conclude that for $0\leq f \in W^{1,2}(\Omega)$ and $\|f\|_{2}=1$ we have
\begin{align}\label{eqn:logsobgen}
\int_\Omega f^2 \log_+f \leq \ep' Q(f)+b_4\beta(\ep')\|f\|_2^2+b_5\int_\Omega |\log d||f|^2.
\end{align}
where $\ep'=b_4\ep>0$ and $\beta(\ep')=b_4\tilde\beta(\ep')$.  Finally
\[
\int_\Omega f^2 \log f \leq \int_\Omega f^2\log_+ f
\]
and given arbitrary $0 \leq f \in W^{1,2}(\Omega)$ with $f\ne 0$ we substitute $f/\|f\|_2$ into (\ref{eqn:logsobgen}) to get the final result.\section{Log--Hardy Inequality}\label{sect:loghardy}

In order to use Theorem \ref{thm:gls} we need to be able to estimate the term
\[
\int_\Omega |\log d(x)||f(x)|^2\dnx.
\]
We do this with a Hardy-type inequality:
\begin{align}\label{eqn:loghardy}
|\log d |^\alpha \leq b_6(H_N+1)
\end{align}
for some constants $\alpha>0$ and $b_6>0$ (this is to be interpreted in the sense of quadratic forms).  We will refer to this as a {\em logarithmic Hardy inequality } (or just a {\em log--Hardy inequality}).

We now assume that $\Omega$ is bounded.  In this case as with the ordinary weak Hardy inequality (see for example \cite{HC}) the log--Hardy inequality depends only on the local geometry of the boundary:

\begin{definition} Let $\Omega$ be a bounded Euclidean domain.  We say that a point $a\in \bd$ is {\em $\alpha$--regular} if there exists a neighbourhood $U$ of $a$ such that
\begin{align}\label{eqn:homer}
\int_\Omega |\log d(x)|^\alpha |f(x)|^2 \dnx \leq \kappa (Q(f)+\|f\|_2^2)
\end{align}
for all $f \in W^{1,2}(\Omega)$ which vanish outside U and some constant $\kappa>0$ which does not depend on $f$.
\end{definition}

\begin{lemma}  If $\Omega$ is bounded and every point of the boundary is $\alpha$-regular then
\[
\int_\Omega |\log d(x)|^\alpha |f(x)|^2 \dnx \leq B(Q(f)+\|f\|_2^2)
\]
for all $f \in W^{1,2}(\Omega)$ and some constant $B>0$ which does not depend on $f$.
\end{lemma}
\Proof  This uses a partition of unity argument.  See for example \cite[Section 2]{HC}.

\begin{remark}
Note that it is not possible in general to prove an inequality of the form
\[
d^{-\gamma} \leq c(H_N+1)
\]
since one could then use the same interpolation argument employed by Burenkov and Davies in \cite{SSNL} to prove an ordinary Sobolev embedding.
\end{remark}

\begin{lemma}\label{lem:eps}
Suppose the log--Hardy inequality (\ref{eqn:loghardy}) holds for some $\alpha >1$.  Then for every $\ep>0$ we have
\[
 \int_\Omega |f(x)|^2|\log d(x)|\dn x \leq \ep Q(f)+ ((\ep/b_6)^{-1/(\alpha-1)}+\ep)\|f\|_2^2
\]
for $f\in W^{1,2}(\Omega)$.
\end{lemma}
\Proof This uses the elementary inequality
\[
t \leq \delta t^{\lambda}+\delta^{-1/(\lambda-1)}
\]
valid for all positive $t>0, \delta>0$ and $\lambda >1$.  Thus
\begin{align*}
 \int_\Omega |f|^2|\log(d(x))|\dn x &\leq \int_\Omega (\delta|\log(d(x))|^{\alpha}+\delta^{-1/(\alpha-1)})|f|^2 \\
& \leq b_6\delta Q(f)+b_6\delta\|f\|_2^2+\delta^{-1/(\alpha-1)}\|f\|_2^2.
\end{align*} Now let $b_6\delta=\ep$. 

\begin{theorem}
Let $\Omega$ be a bounded region and suppose the log--Hardy inequality (\ref{eqn:loghardy}) holds for some $\alpha  >1$.

  Then we have the log--Sobolev inequality
\[
\int_\Omega f^2\log f \dn x \leq \ep Q(f)+\eta (\ep)\|f\|^2_2+\|f\|_2^2\log \| f\|_2
\]
for all $f\in W^{1,2}(\Omega)$ where $0<\ep$ and
\[
\eta(\ep)\leq b_7\ep^{-1/(\alpha\beta -1)}-b_8\log \ep +b_9.
\]
for some constants $b_7,b_8$ and $b_9>0$
\end{theorem}

\Proof

By Theorem \ref{thm:gls} we have
\[
 \int_{\Omega}f^2\log f \leq \ep Q(f)+\beta(\ep)\|f\|^2_2+\|f\|_2^2\log \| f\|_2+b_0\int_{\Omega}|\log d| f^2.
\]

Applying Lemma \ref{lem:eps} we have
\begin{align*}
\int_\Omega f^2\log f \dnx &\leq \ep'Q(f)+\beta(\ep')\|f\|_2^2+\|f\|_2^2\log \|f\|_2 \\
&+b_0\ep'Q(f)+b_0((\ep'/b_6)^{-1/(\alpha\beta-1)}+\ep')\|f\|_2^2.
\end{align*}

Now let 
\[
\eta(\ep')=\beta(\ep')+b_0\ep'+b_0b_6(\ep')^{-1/(\alpha\beta-1)}
\]
and the result follows by scaling $\ep'$.\section{Euclidean Domains with Exponential Cusps}\label{sect:euclideandomains}

We now give an application of the previous results.  The domains that we consider will be simple in order to make the general method clear but further applications are possible.\begin{definition}
Let $\Omega \subset \rn$ be a bounded domain.  Then we say that $\bd$ is $(\log,\alpha)$ {\em regular} or more tediously it  has a logarithmic modulus of continuity with exponent $\alpha$ near $a\in \bd$
  if $a$ has a neighbourhood $U$ which can be represented in the following form after translation and rotation of coordinates.  The set $U$ is of the form
\[
U:=\{ (x',x_N):   x' \in B \text{ and } 0< x_N < g(x')\}
\]
where $B$ denotes the ball
\[
B:=\{ x \in \mathbb{R}^{N-1} : |x| < 1/2\}
\]
and  $0<g$ is a function that satisfies
\begin{align}\label{prop:f}
|g(x')-g(y')| \leq A | \log(|x'-y'|)|^{-\alpha}
\end{align}
for all $x',y'\in B$ and $\alpha>0$.
\end{definition}
We also define $\Gamma$ to be the set
\[
\Gamma := \{ (x',g(x') : x'\in B\}
\]
and let the function $d_\Gamma$ be defined by
\[
d_\Gamma(x)=\text{dist}(x,\Gamma):=\inf_{y \in \Gamma}|x-y|
\]
for all $x \in \Omega$.  Note that $d_\Gamma(x)=d(x)$ for all $x$ in a sufficiently small neighbourhood of $a$.

The following lemma is a modification of the case when the boundary function $g$ is assumed to be H\"older continuous..

\begin{lemma}\label{lem:ed}

Let $\bd$ be $(\log,\alpha)$ regular near $a$.  For $x=(x',x_N)\in U$ define the function $e(x)$ by
\[
e(x)=g(x')-x_N.
\]
Then for all $x \in U$ such that
\[
e(x)<(1+A)^{-\alpha}
\]
we have
\[
\exp\left(-\frac{e(x)^{-1/\alpha}}{1+A}\right)\leq d_\Gamma(x) \leq e(x).
\]
\end{lemma}

\Proof This proof follows the H\"older case in \cite[Lemma 4.6]{K}.
Given $z=(z',z_N)\in \bd$ we define the cusp $C(z)$ by
\[
C(z):=\{ (x',x) | x' \in B \text{ and } x_N < g(z')-A|\log|z'-x'||^{-\alpha} \}.
\]

Now let $(x',x)\in C(z)$.  From the property (\ref{prop:f}) we have
\[
g(x') \geq g(z')-A|\log |z'-x'||^{-\alpha}
\]
and thus
\[
0<x_N< g(z')-A|\log |z'-x'||^{-\alpha} \leq g(x')
\]
which implies that $(x',x_N) \in U$.  Thus we have shown that
\begin{align}\label{eqn:subset}
C(z) \subset U.
\end{align}
Now, given $x=(x',x_N)\in \Omega$ we define the constant $R$ by
\[
R:=\exp\left(-\frac{e(x)^{-1/\alpha}}{1+A}\right)
\]
and then consider the closed ball $B(x,R)$ given by
\[
B(x,R):=\{U \in \Omega:|x-y|\leq R\}.
\]
Thus, for $y \in B(x,R)$ we have $|x'-y'| \leq R$ and $|x_N-y_N|\leq |\log R|^{-\alpha}$.

Now
\begin{align*}
y_N-g(x')+A|\log |x'-y'||^{-\alpha}&=y_N-x_N+x_N -g(x')+A|\log |x'-y'||^{-\alpha} \\
&\leq |y_N-x_N|-(g(x')-x_N)+A|\log |x'-y'||^{-\alpha} \\
&\leq (1+A)|\log R|^{-\alpha}-e(x) \\
&=0
\end{align*}
Hence $B(x,R) \subset \overline{C(z)}$.  Combining this with (\ref{eqn:subset}) and the result follows.

\begin{theorem}\label{thm:ws}
Let $\bd$ be $(\alpha,\log)$ regular near $a$.

Let $f \in W^{1,2}(\Omega)$ and $f$ vanish outside a neighbourhood of $a$.  Then there exists $b_{10}<\infty$ such that
\[
\int_\Omega |f|^2|\log(d_\Gamma(x))|^{\alpha\beta}\dn x \leq b_{10}
  \left( \int_\Omega (|\nabla f|^2 + |f|^2)\dn x \right)
\]
for all $0<\beta<1$.
\end{theorem}

\Proof If $0<\beta<1$ then the embedding $W^{1,2}(I)\subseteq L^\infty(I)$ for any finite interval $I$ implies that there exists $b_{10}>0$ such that
\[
\int_\Omega e(x)^{-\beta} |f(x)|^2 \dn x \leq b_{10} \int_\Omega (|\nabla f|^2+|f|^2)\dn x
\]
for functions supported in a neighbourhood of $\Gamma$.  Now apply Lemma \ref{lem:ed} and the result follows immediately.

\begin{corollary}  Let $\bd$ be $(\log,\alpha)$ regular near every point and let $0<\beta<1$ be such that $\alpha\beta>1$.  Then there exists $b_{11}>0$ such that the log--Hardy inequality
\[
|\log d|^{\alpha\beta} \leq b_{11}(H_N+1)
\]
holds in the sense of quadratic forms.  Moreover we have the log--Sobolev inequality  \[
\int_\Omega f^2\log f \dn x \leq \ep Q(f)+\eta (\ep)\|f\|^2_2+\|f\|_2^2\log \| f\|_2
\]
for all $f\in W^{1,2}(\Omega)$ where $0<\ep$ and
\[
\eta(\ep)\leq b_{12}\ep^{-1/(\alpha\beta -1)}-b_{13}\log \ep +b_{14}.
\]
for some constants $b_{12},b_{13}$ and $b_{14}>0$.
\end{corollary}

\section{A Rotationally Invariant Example}\label{sect:rot}

The previous results cover the case in which the cusp is not too sharp (i.e. $\alpha>1$).  It would be interesting to know what happened outside this range of $\alpha$. 

The following section considers a slightly different problem.  We  now work on a two dimensional, rotationally invariant Riemannian manifold.  The manifold when embedded into $\mathbb{R}^3$ has a cusp either at a finite point or at infinity.  This model was considered by Davies in \cite{HKRM} .  Davies showed ultracontractivity in the case of Dirichlet boundary conditions with a similar function $\eta(\ep)$ to that which we found in the previous section.  However, at some critical value of the parameter controlling the sharpness of the cusp we show that ultracontractivity fails by demonstrating that one of the eigenfunctions does not lie in $L^\infty$ (this of course in no way contradicts the compactness of the resolvent).  This is possible because the rotational invariance of the problem allows us to reduce it to one that is one-dimensional.

This breakdown in ultracontractivity is interesting in its own right and although it does not allow us to deduce anything about the flat case it does show that ultracontractivity can break down before we have exhausted all possible functions $\eta(\ep)$.

\subsection{Basic Model}
We begin by recalling the definition of the manifold from \cite[Example 15]{HKRM} and a few elementary facts about it.

 Let $M$ be the manifold
\[
M:=(2\pi, \infty) \times S^1
\]
equipped with the metric
\[
ds^2=g(u)(\id u^2+\id \theta^2).
\]
Thus, the Riemannian volume element $\dvol$ is given by
\[
\dvol = g(u)\id u \id \theta.
\]
If $g$ is bounded and
\[
\int_{2\pi}^\infty g(u)^{1/2}\id u <\infty
\]
then $M$ is bounded and if
\[
|g'|<2g
\]
for large $u$, then we may embed the manifold in $\mathbb{R}^3$ for such $u$ by setting
\[
\left\{\begin{array}{c}x=g(u)^{1/2}\cos \theta \\ y=g(u)^{1/2} \sin \theta \\ z=z(u) \end{array} \right.
\]
where
\[
z'=\sqrt{g-(g')^2/4g}.
\]

One gets a power cusp by setting $g(u)=u^{-\alpha}$ where $\alpha>2$.  The next case Davies introduces involves setting
\begin{align}\label{eqn:basicmetric}
g(u)=u^{-2}(\log u)^{-\alpha}.
\end{align}
The manifold $M$  has finite volume for all $\alpha>0$,     however, it is bounded if and only if $\alpha>2$.

For all $\alpha>0$ the curvature $K$ (see \cite[Example 15]{HKRM}) has the asymptotic behaviour
\[
K\sim -(\log u)^\alpha \text{ as } u \to \infty.
\]
\begin{lemma}
Let $\alpha>2$ and let $B$ be the ball centred at the cusp with radius $\ep$.  Let $\kappa:=(\alpha-2)/2$. 
  Then we have
\[
\text{Vol}(B)\sim \exp \left(- \kappa^{-1/\kappa} \ep^{-1/\kappa} \right)(\kappa\ep)^{-1/\kappa}
\]
as $\ep \to 0^+$.
\end{lemma}
\Proof  The distance between any point $(u_0,\theta_0)$ and the cusp is now
\[
\int_{u_0}^\infty u^{-1}(\log u)^{-\alpha/2}\id u = \frac{2}{\alpha-2}(\log u_0)^{1-\alpha/2}.
\] 
To simplify notation let
\[
\eta(\ep):=\exp \left( \kappa^{-1/\kappa} \ep^{-1/\kappa}\right).
\]
The volume of the ball is then
\begin{align*}
\int_{\eta(\ep)}^\infty u^{-2}(\log u)^{-\alpha} \id u &=\int_{\log \eta(\ep)}^\infty e^{-v}v^{-\alpha} \id v \\
&\sim \frac{(\log \eta(\ep))^{-\alpha}}{\eta(\ep)} \text{ as } \ep \to 0^+.
\end{align*}

\subsection{The Proof of Ultracontractivity}

The main result in \cite[Example 15, B]{HKRM} is the following:

\begin{theorem}\cite[Example 15, Case B]{HKRM}\label{thm:mainsemi}  Let the metric be given by (\ref{eqn:basicmetric}) and suppose $\alpha>2$.  Then the manifold is bounded and we have a logarithimic Sobolev inequality
\[
\int_M f^2\log f \dvol \leq \ep Q(f)+\beta(\ep)\|f\|_2^2+\|f\|_2^2\log \|f\|_2
\]
for all $0\leq f \in W^{1,2}_0(M)$ and $0<\ep<1$ where
\[
\beta(\ep)\leq c\ep^{-1/(\alpha-1)}.
\]
\end{theorem}

Davies' technique involves two components.  The first is the existence of the quadratic form inequality
\begin{align}\label{eqn:keyquad}
(\log u)^\alpha \leq c_0H_D
\end{align}
where $H_D$ is the Dirichlet Laplacian.  The second component involves finding a uniform covering of the manifold with sets $\Omega_n$ that are diffeomorphic to cubes and also uniform estimates for the metric and other quantities in each $\Omega_n$.  This procedure is explained in \cite[Section 3]{HKRM}.

However, if one studies the proof of Theorem \ref{thm:mainsemi} the fact that the manifold is bounded is unimportant.  Moreover, the quadratic form inequality (\ref{eqn:keyquad}) holds for all $\alpha>0$. Thus we have:

\begin{theorem}\label{thm:alphalessone}
 Let the metric be given by (\ref{eqn:basicmetric}) and suppose $2\geq\alpha>1$.  Then the manifold is unbounded but of finite volume and we have a logarithmic Sobolev inequality
\[
\int_M f^2\log f \dvol \leq \ep Q(f)+\beta(\ep)\|f\|_2^2+\|f\|_2^2\log \|f\|_2
\]
for all $0\leq f \in W^{1,2}_0(M)$ and $0<\ep<1$ where
\[
\beta(\ep)\leq c\ep^{-1/(\alpha-1)}.
\]
\end{theorem}

In fact we actually have more than this because (\ref{eqn:keyquad}) actually holds with $H_D$ replaced by $H_N$.  
\begin{lemma}  Let ${\cal D}$ denote the set
\[
{\cal D}:=\{ f \in C^\infty([2\pi,\infty)\times S^1)\cap W^{1,2}(M): \text{supp}(f)\cap [R_f,\infty)=\emptyset \text{ for some } R_f>0 \}.
\]
Then ${\cal D}$ is dense in $W^{1,2}(M)$.
\end{lemma}

\Proof We have that $W^{1,2}(M,u\id u \id \theta)\hookrightarrow W^{1,2}(M,\dvol)$ and since $M$ has the segment property \cite[Theorem 3.18]{A} gives us that
\[
W^{1,2}(M, u\id u \id \theta  )=\overline{{\cal D}}.
\]
We  now prove a version of the result of Moss, Allegretto and Piepenbrink used to prove (\ref{eqn:keyquad}), suitable for the case of Neumann boundary conditions.\begin{theorem}
We have
\begin{align}\label{eqn:lisa}
H_N \geq \frac{3}{16}(\log u)^\alpha
\end{align}
in the sense of quadratic forms.
\end{theorem}

\Proof  We define $\phi$ by
\[
\phi(u):=u^{1/2}-Cu^{1/4}
\]
where $C:=2(2\pi)^{1/4}$.  Thus $\phi$ satisfies the Neumann condition $\phi'(2\pi)=0$.  Then for $f\in {\cal D}$ we have
\begin{align*}
\int_{2\pi}^\infty \frac{\partial \phi}{\partial u}\frac{\partial f}{\partial u} \id u &=\int_{2\pi}^\infty \left(\hf u^{-1/2}-\frac{C}{4}u^{-3/4}\right)\frac{\partial f}{\partial u} \id u \\
=&\left. f\left(\hf u^{-1/2}-\frac{C}{4}u^{-3/4}\right)\right|^\infty_{2\pi}-\int_{2\pi}^\infty f\left(-\qt u^{-3/2}+\frac{3C}{16}u^{-7/4}\right) \id u\\
&=\int_{2\pi}^\infty f\left(\qt u^{-3/2}-\frac{3C}{16}u^{-7/4}\right) \id u \\
&=\int_{2\pi}^\infty(\log u)^\alpha f\left(\qt u^{1/2}-\frac{3C}{16}u^{1/4}\right)u^{-2}(\log u)^{-\alpha} \id u
\end{align*}
and so
\begin{align*}
\int_M            \left( \frac{\partial \phi}{\partial u}\frac{\partial f}{\partial u}+\frac{\partial \phi}{\partial \theta}\frac{\partial f}{\partial \theta}\right) \id u\id \theta
&\geq \frac{3}{16} \int_M(\log u)^\alpha f\phi\dvol.
\end{align*}

We now follow the proof of \cite[Theorem 1.5.12]{HKST}.  Given $f\in {\cal D}$ we set $f=\phi g$ for $g\in {\cal D}$.  Then
\begin{align*}
\int_M \left( \left|\frac{\partial f}{\partial u}\right|^2+ \left|\frac{\partial f}{\partial \theta}\right|^2\right) \id u \id \theta &\geq \int_M \left( |g|^2\left|\frac{\partial \phi}{\partial u}\right|^2+2g \frac{\partial \phi}{\partial u} \right) \id u \id \theta \\
&=\int_M \frac{\partial \phi}{\partial u}\frac{\partial \phi|g|^2}{\partial u}\id u\id \theta \\
&\geq \frac{3}{16}\int_M (\log u)^\alpha \phi^2 |g|^2 \dvol \\
&=\frac{3}{16}\int_M (\log u)^\alpha |f|^2 \dvol.
\end{align*}

Thus we have
\[
H_N \geq \frac{3}{16}(\log u)^\alpha
\]
in the sense of quadratic forms.
 
\begin{corollary}  Theorems \ref{thm:mainsemi} and \ref{thm:alphalessone} hold also in the case of Neumann boundary conditions.
\end{corollary}

\Proof  The domain decomposition does not depend on the particular boundary conditions chosen.  Its critical use is in \cite[Theorem 8]{HKRM}.  However,  the proof of this depends on the classical bound
\[
h \leq \alpha_1\|h\|_{2N}(K+1)^{1/4}
\]
where $\alpha_1$ is some constant, $h \geq 0$ and $K$ is the {\em Neumann} Laplacian on the cube $(0,1)^N$ (all of this is to be interpreted in the Euclidean metric).  It is easily seen that all of the calculations hold not only for $f\in C^\infty_c(M)$ but also for $f\in \cal{D}$.  Since we have both the conclusion of \cite[Theorem 8]{HKRM} and the quadratic form inequality (\ref{eqn:lisa}) we can apply \cite[Theorem 2]{HKRM} which is stated for an abstract positive self--adjoint operator to get the final result.

\subsection{Breakdown of Ultracontractivity}

Our main result shows that ultracontractivity breaks down at the point $\alpha=1$.
\begin{theorem}
Let the metric be given by (\ref{eqn:basicmetric}) and let $\alpha=1$.  Let $H$ be the Laplacian subject either to Dirichlet or Neumann boundary conditions.  Then $e^{-Ht}$ is {\em not} ultracontractive.
\end{theorem}

\Proof  Suppose that $e^{-Ht}$ {\em is} ultracontractive.  Then since the volume of $M$ is finite we deduce that $H$ has compact resolvent.  Now, since the domain is rotationally invariant so we may use the rotational group to decompose $L^2(M)$ into orthogonal linear spaces $\{L_n\}_{n\in\mathbb{Z}}$ consisting of functions of the form
\[
h(r\cos \theta, r \sin \theta)=f(r)e^{in\theta}.
\] 
Since the operator commutes with rotations it maps each of these subspaces into itself and so its spectral behaviour can be analysed in each subspace independently.  This allows us to reduce $H$ to a one dimensional operator.  The associated differential equation is then

\begin{gather}
-f''(u)+n^2f(u)=\lambda g(u)f(u), \label{eqn:model1}   \\
2\pi < u < \infty, \notag \\
n \in \mathbb{Z}. \notag
\end{gather}

The appropriate boundary conditions are given by classifying the end points.  The left hand end point $u=2\pi$ is clearly regular and we may specify boundary conditions in the usual way.  Now let $\lambda=0$.  Then the equation has the basis of solutions
\[
\begin{array}{cc} \phi_1(u)= 1 & \phi_2(u)=u
\end{array}
\]
in the case $n=0$.  If $n\ne 0$ then we have the following basis of solutions:
\[
\begin{array}{cc} \psi_1(u)= e^{-nu} & \psi_2(u)=e^{nu}
\end{array}
\]

A calculation shows that 
\[
\phi_2,\psi_2 \not \in L^2((2\pi,\infty),g(u)\id u)
\]
 whereas 
\[
\phi_1, \psi_1 \in L^2((2\pi,\infty),g(u)\id u)\]
whence we classify $\infty$ as being Limit Point (for an introduction to the theory of singular Sturm-Liouville problems and end point classification see \cite{NZ92}).

Let
\begin{align*}
M(2\pi,\infty):=\{f: f,f'\in AC(2\pi,\infty) \text{ and } f,f'' \in L^2(2\pi,\infty, g(u),\id u) \}
\end{align*} 
and then the domain of the one dimensional operator subject to Dirichlet boundary conditions is
\[
{\cal D}_D=\{f \in M(2\pi,\infty) : f(2\pi)=0\}
\]
and subject to Neumann conditions:
\[
{\cal D}_N=\{f \in M(2\pi,\infty) : f'(2\pi)=0\}.
\]
We now focus on the subspace of purely radial functions by taking  $n=0$ and make the change of variable $\log u=v$.  Let $h(v)=f(1/u)$ and the equation becomes
\begin{gather*}
-(e^{-v}h'(v))'=\frac{\lambda e^{-v} h(v)}{v} \\
\log 2\pi < v < \infty.
\end{gather*}
Now let $k(v)=e^{-v/2}h(v)$ and the equation now becomes
\begin{gather*}
-k''(v)+V(v)k(v)=0 \\
\log 2\pi < v < \infty
\end{gather*}
where
\[
V(v)=\qt-\frac{\lambda}{v}.
\]
We can now use standard techniques from asymptotic analysis (see e.g. \cite[Chapter 6]{M}) to analyse this equation.  The solutions $k_1$ and $k_2$  have the following asymptotic forms:
\begin{align*}
k_1(v)=Ae^{-v/2}v^\lambda(1+O(1/v)) \\
k_2(v)=Be^{v/2}v^{-\lambda}(1+O(1/v))
\end{align*}
as $v \to \infty$ and where $A$ and $B$ are constants.  Consequently the solutions to equation (\ref{eqn:model1}) satisfy
\begin{align*}
f_1(u)=Au(\log u)^{-\lambda}(1+O(1/(\log(u))) \\
f_2(u)=B(\log u)^{\lambda}(1+O(1/(\log(u)))
\end{align*}
as $u\to \infty$.

Thus $f_1 \not \in L^2((2\pi,\infty),g(u)\id u)$ and $f_2  \in L^2((2\pi,\infty),g(u)\id u)$ but also neither $f_1$ nor $f_2$ belong to $L^\infty$ (unless of course $\lambda=0$).  However, the assumption that $e^{-Ht}$ is ultracontractive implies that all eigenfunctions lie in $L^\infty$.  Thus we either have unbounded eigenfunctions or no eigenfunctions both of which contradict the ultracontractivity assumption.

\textbf{Acknowledgements}  I would like to thank Brian Davies for suggesting this problem and for his guidance and advice during this work.  I also acknowledge the support of the Engineering and Physical Sciences Research Council through a research studentship.

 Department of Mathematics \\
King's College London \\
Strand \\
London \\
WC2R 2LS \\
U.K. \\
cmason@mth.kcl.ac.uk
\end{document}